\documentclass{article}

\usepackage{amsmath}
\usepackage{url}
\RequirePackage{graphicx}
\RequirePackage{pifont,latexsym,ifthen,amsthm,rotating,calc,booktabs}
\RequirePackage{amsfonts,amssymb,amsbsy,amsmath}

\theoremstyle{definition}
\newtheorem{theorem}{Theorem}

\newtheorem{corollary}[theorem]{Corollary}
\newtheorem{proposition}[theorem]{Proposition}
\newtheorem{definition}[theorem]{Definition}
\newtheorem{example}[theorem]{Example}

\newcommand\lb{[\![}
\newcommand\rb{]\!]}

\makeatletter
\@twosidetrue
\makeatother

\begin{document}

\author{Jonathan Fine\\
15 Hanmer Road, Milton Keynes, MK6 3AY, United Kingdom\\
email: jfine@pytex.org}
\title{A complete $h$-vector for convex polytopes}
\date{30 November 2009}

\maketitle

\begin{abstract}
\noindent
This note defines a complete $h$-vector for convex polytopes, which
extends the already known toric (or mpih) $h$-vector and has many
similar properties.  Complete means that it encodes the whole of the
flag vector.

First we define the concept of a generalised $h$-vector and state some
properties that follow.  The toric $h$-vector is given as an example.
We then define a complete generalised $h$-vector, and again state
properties.  Finally, we show that this complete $h$-vector and all
with similar properties will sometimes have negative coefficients.

Most of the proofs, and further investigations, will appear elsewhere.
\end{abstract}

\section{Generalised $h$-vectors}

This note defines a complete $h$-vector for convex polytopes, and
states some of its properties.  Background, motivation and most of the
proofs will be given elsewhere~\cite{fine:chv2}.  Prior knowledge of
the toric $h$-vector, for example as in \cite{sta87:gen} or
\cite[\S4.1]{de2009decomposition}, will help the reader.  Throughout
$\Delta$ will denote a convex polytope of dimension $d$.  We study
linear functions $h=h(\Delta)$ of the flag vector $f = f(\Delta)$ of
$\Delta$.  When conversely $f(\Delta)$ can be computed from
$h(\Delta)$ we say that $h$ is \emph{complete}.

Let $\delta$ be a face of $\Delta$, of dimension~$i$.  Associated with
$\delta\subseteq\Delta$ there is the \emph{link} $L_\delta$, a convex
polytope of dimension $d - i -1$, which encodes the local geometry of
$\Delta$ around $\delta$.  For example, if $\delta$ is a vertex then
around that vertex $\Delta$ looks like $CL_\delta$, when $C$ is the
\emph{cone} or pyramid operator.  Although $L_\delta$ is determined
only up to projective equivalence, its flag vector is an invariant of
$\delta\subseteq\Delta$.  It is convenient to set $L_\Delta =
\emptyset$ and $C\emptyset = \mathrm{pt}$.

Throughout $C$ and $I$ denote the cone and cylinder (or pyramid and
prism) operators, and we think of
\begin{equation}
  D = IC - CC
\end{equation}
and $I$ and $C$ as operators on flag vectors. The \emph{total link
  vector} $\ell = \ell(\Delta)$ has components $\ell_i =\sum_{\dim
  L_\delta=i} f(L_\delta)$.  Many of the results rely on
$\ell(I\Delta) = (1 + 2C)\ell(\Delta)$ and $\ell(C\Delta) = (1 +
C)\ell(\Delta) + f(\Delta)$.  In~\cite{fine:mv+ic} the author shows
$DI=ID$.  This is a partial expression of the next result, upon which
the definition of $h$ relies.

\begin{theorem}[{Bayer-Biller~\cite{baybil85:gen}, generalised Dehn-Sommerville}]
\label{thm:gds}
The flag vectors produced by applying all words $W$ in $C$ and $D$ to
a point are a basis for the vector space spanned by all convex
polytope flag vectors.
\end{theorem}

\begin{proposition}
\label{prop:linear}
Suppose $g$ is a linear function of flag vectors.  Then
\begin{equation}
  h(\Delta) = \sum\nolimits_{\delta\subseteq\Delta}(x-y)^{\dim\delta}g(L_\delta)
\end{equation}
is also a linear function of flag vectors.
\end{proposition}

\begin{definition}
If $h(\Delta)$ is as in Proposition~\ref{prop:linear} and in addition
both of
\begin{enumerate}
\item $g(\emptyset) = 1$ (or equivalently $h(\mathrm{pt}) = 1$).
\item $g(CL) = y\,g(L)$.
\end{enumerate}
hold then we will say that $h$ is a \emph{generalised $h$-vector}.
\end{definition}

\newpage
\begin{proposition}
\label{prop:ghv}
Suppose $h$ is a generalised $h$-vector and $\Delta$ is a convex
polytope.  Then
\begin{enumerate}

\item $h(I\Delta) = (x + y)\, h(\Delta)$.

\item $h(C\Delta) = g(\Delta) + x\,h(\Delta)$, and thus we can define
  $g$ from $h$.

\item $h(D\Delta) = xy\,h(\Delta)$ (which follows from the two
  previous statements).

\item If $\Delta$ is simple then $h(\Delta) = \sum_{i=0}^d
  (x-y)^iy^i\,f_i(\Delta)$, the usual formula for simple polytopes.

\item If $\Delta_1$ is simple and $\Delta_2$ is any convex polytope
  then $h(\Delta_1\times\Delta_2) = h(\Delta_1)h(\Delta_2)$.

\end{enumerate}
\end{proposition}

\begin{corollary}
Suppose there is a complete generalised $h$-vector.  Then as operators
on flag vectors $DI =ID$.
\end{corollary}

\begin{proof}
We have $h(DI\Delta) = h(ID\Delta)$.  If $h$
is complete $f$ is a linear function of $h$ and so $f(DI\Delta) =
f(ID\Delta)$.
\end{proof}

\begin{proposition}
Suppose $h$ is a linear function of flag vectors with
$h(\mathrm{pt})=1$, $h(I\Delta) = (x+y)\,h(\Delta)$ and $h(D\Delta) =
xy\,h(\Delta)$ for any convex polytope $\Delta$.  Then $h$ is a
generalised $h$-vector (with $g(\Delta) = h(C\Delta) -xh(\Delta)$).
\end{proposition}

\begin{proposition}
\label{prop:gDW}
Suppose $g(DW\mathrm{pt})$ is known for all words $W$ in $C$ and $D$.
Then $g$ determines a unique generalised $h$-vector. (This follows
from Theorem~\ref{thm:gds}, $g(CL)=y\,g(L)$, $g(\emptyset)=1$ and
linearity.)
\end{proposition}

\begin{proposition}
The formula
\begin{equation}
\label{eqn:mpihg}
  g(DL) = xy\,g(L)
\end{equation}
defines the toric (or middle perversity intersection homology or mpih)
$h$-vector, as in \cite{sta87:gen} or
\cite[\S4.1]{de2009decomposition}.
\end{proposition}

\section{A complete $h$-vector}

For simple polytopes $h(x,y) = h(y,x)$ or in other words $h$ is
\emph{palindromic}. This is a very important property.  We will use
the following notation.  We denote, for example, $ay^2 + bxy +c^2x$ by
$[a,b,c]$, and $[1,1,1]$ by $\lb 2\rb$.  We denote $(xy)^i\lb j\rb$ by
$\lb i,j\rb$.  Thus $\lb 0,j\rb = \lb j\rb$ and $\lb1,0\rb = [0,1,0] =
xy$.

\begin{definition}[Keyed and palindromic generalised $h$-vectors]

Suppose a generalised $h$-vector has the form
\begin{equation}
  h(\Delta) = \sum h_k(\Delta)w_k
\end{equation}
where each $h_k$ is a homogeneous polynomial, each $k$ is a \emph{key}
as defined below and $w_k$ its associated symbol.  Suppose also that
$\dim\Delta = \deg h_k + \deg k$.  If all this holds, we say that $h$
is a \emph{keyed} generalised $h$-vector.  If each $h_k$ is
palindromic we say that $h$ is \emph{palindromic}.
\end{definition}

For the rest of this note $h$ denotes the palindromic keyed
generalised $h$-vector we are about to define.  Recall that by
Proposition~\ref{prop:gDW} it is enough to define $g$ on polytopes of
the form $g(DW\mathrm{pt}$).

\begin{definition}[Complete generalised $h$-vector]
\label{def:g}
Suppose $h(v) = \lb i,j\rb w_k$.  Then
\begin{equation}
  g(Dv) = (xy)^{i+1}\,y^{j+1}\, w_k \>+ \>w_{k'}
\end{equation}
where $k'$ is as below.  This we extend linearly to
$h(v)=\sum\lambda_{ijk}\lb i,j\rb w_k$ and so to $v = W\mathrm{pt}$.
\end{definition}

It is clear that $\deg k' = \deg k + 2i + j + 3$.  The definition
assumes that $h(W\mathrm{pt})$ is palindromic.  It turns out to be the
same as (\ref{eqn:mpihg}), except for the addition of $w_{k'}$.  This
addition has, of course, recursive consequences.  The author has
developed software \cite{fine:python-hvector} for $h$-vector
computations.  In particular \cite{fine:python-hvector} contains a
human and machine readable table of $h$-vectors for $CD$ polytopes up
to dimension~$10$.

\begin{definition}[$h$-key]
If $k = ((d_1, \ldots, d_r), (c_1, \ldots, c_r))$ then $k' = ((i, d_1,
\ldots, d_r), (j, c_1, \ldots, c_r))$.  As a shorthand we sometimes
write, for example, $k=((1,3,2),(0,2,1))$ as $132{;}021$.  We write
$e$ for the empty key $((),())$, and set $w_e =1$.  Thus,
$h(\mathrm{pt}) = \lb0,0\rb w_e = 1$.  We use $\deg k = 2\sum d_i +
\sum c_i + 3r$ to define the \emph{degree} of $k$.
\end{definition}

\begin{example}
Suppose $h(v) = [1,1]w_k$, which we can write as $\lb 0,1 \rb w_k$.
Then
\begin{align}
  g(Dv) &= (xy)^1y^2w_k + w_{k'} = [0,1,0,0,0]w_k + w_{k'}\\
  h(CDv) &= g(DV) + x h(Dv) = [0,1,0,0,0]w_k + w_{k'} + [0,0,1,1,0]w_k\\
  & = [0,1,1,1,0]w_k + w_{k'} = \lb 1,2\rb w_k + w_{k'}\>.
\end{align}
\end{example}

\begin{proposition}
If $h(v) = \lb i,j\rb w_k$ then $h(CDv) = \lb i+1, j+1\rb w_k + w_{k'}$.
\end{proposition}

\begin{proposition}
\label{prop:hcv}
If $h(Cv) = \sum \lambda_{ijk}\lb i,j\rb w_k$ then $h(CCv) = \sum
\lambda_{ijk}\lb i,j+1\rb w_k$.
\end{proposition}

\begin{example}
If $h(v)=w_{0;0}$ then $v=(CD-DC)\mathrm{pt}$, because
$h(CD\mathrm{pt}) = \lb 1,1\rb + w_{0;0}$ and $h(DC\mathrm{pt})=\lb
1,1\rb$.  We write $w_{0;0}$ as $\lb 0,0\rb w_{0;0}$.  We now have
\begin{align}
h(CCD\mathrm{pt}) &= \lb 1,2 \rb + \lb 0, 1 \rb w_{0;0}\\
h(CDC\mathrm{pt}) &= \lb 1,2 \rb + w_{0;1}\\
h(Cv) &= \lb 0, 1 \rb w_{0;0} - w_{0;1} \>.
\end{align}
This shows that sometimes $h(v)=\lb i,j\rb w_k$ does not imply $h(Cv) = \lb
i, j+1\rb w_k$.
\end{example}

\begin{proposition}[Converse to Proposition~\ref{prop:hcv}]
If $h(Cv) = \sum \lambda_{ijk}\lb i,j+1\rb w_k$ then $v = Cv'$ for
some $v'$.
\end{proposition}

\begin{theorem}
\label{thm:ghv}
Let $g$ be as in Definition~\ref{def:g}.  Then
\begin{enumerate}

\item There is a unique extension of $g$ to a generalised $h$-vector
  $h$.

\item $h(\Delta)$ is palindromic.

\item $h_e$ is the toric $h$-vector.

\item The matrix for $g$, for the $CD$ and $\lb i,j\rb w_k$ bases, is
  upper triangular with ones along the diagonal.

\item $h$ is complete.

\item Theorem~\ref{thm:gds} can be proved as part of the $g$-$h$
  recursion.

\end{enumerate}
\end{theorem}

Unlike the toric/mpih $h$-vector, the complete $h$-vector $h$ can have
negative coefficients. This is unavoidable.

\begin{example}[Bayer, personal communication]
\label{example:bayer}
Let $P$ be the bipyramid on the $3$-simplex.  Then
\begin{equation}
  h(P) = [1, 4, 10, 4, 1] + [6, 6] w_{0;0} + [-4]w_{0;1}
\end{equation}
\end{example}

\begin{proposition}
Suppose $h'$ is a complete palindromic keyed generalised $h$-vector.
Let $P$ be as in Example~\ref{example:bayer} and let $Q$ be
$CICC\mathrm{pt}$.  Write $h'_{0;1}(P) = [a]$ and $h'_{0;1}(Q) = [b]$.
Then $a>0$ and $b<0$ or vice versa.
\end{proposition}

\begin{proof}[Sketch of proof]
The general form of $h'(\Delta)$ for $\dim\Delta=3$ is
$[a_0,a_1,a_1,a_0] + [b_0]w_{0{;}0}$ and Proposition~\ref{prop:ghv}
determines $h'$ on $CCC\mathrm{pt}$ and $DC\mathrm{pt}$.  Both have
$b_0=0$ and so, for $h'$ to be complete, $h'(CD\mathrm{pt})$ has
$b_0\ne 0$.

Similarly the general form of $h'(\Delta)$ for $\dim\Delta=4$ is
$[a_0,a_1,a_2,a_1,a_0] + [b_0, b_1]w_{0{;}0} + [c_0]w_{0{;}1} $ and
Proposition~\ref{prop:ghv} determines $h'$ on $CCCC\mathrm{pt}$
$DCC\mathrm{pt}$, $DD\mathrm{pt}$.  They all have $c_0=0$.
Proposition~\ref{prop:hcv} shows the same for $CCD\mathrm{pt}$.

For $h'$ to be complete $c_0$ must be non-zero for some convex
polytope.  Think of $c_0$ as a linear function on $4$-polytope flag
vectors.  It is non-zero and vanishes on the hyperplane $H$ spanned by
$f(CCCC\mathrm{pt})$ $f(DCC\mathrm{pt})$, $f(DD\mathrm{pt})$ and
$f(CCD\mathrm{pt})$.  Finally, by a calculation we omit, $f(P)$ and
$f(Q)$ are separated by $H$.  The result follows.
\end{proof}

\bibliographystyle{amsplain}
\bibliography{ehv-note2.bib}

\end{document}